\documentclass[12p]{amsart}
\usepackage{amssymb}
\usepackage{amsmath}
\usepackage{amsfonts}
\usepackage{geometry}
\usepackage{graphicx}
\usepackage{mathrsfs,amssymb}

\usepackage{hyperref}
\usepackage{cleveref}

\theoremstyle{plain}

\newtheorem{lemma}{Lemma}

\newtheorem{proposition}{Proposition}

\newtheorem{theorem}{Theorem}
\numberwithin{equation}{section}

\begin{document}

\title[]{The Nonlinear Schr\"odinger equation on Z and R with bounded initial data: examples and conjectures}

\author{Benjamin Dodson}
\address{Department of Mathematics, Johns Hopkins University, 3400 N. Charles St., Baltimore, MD, 21218}
\email{bdodson4@jhu.edu}
\author{Avraham Soffer}
\address{Department of Mathematics, Rutgers University, 110 Frelinghuysen Rd., Piscataway, NJ 08854}
\email{soffer@math.rutgers.edu}
\author{Thomas Spencer}
\address{School of Mathematics, Institute for Advanced Study, Princeton, NJ, 08540}
\email{spencer@ias.edu}

\begin{abstract}
We study the nonlinear Schr\"odinger equation (NLS) with bounded initial data which does not vanish at infinity.  Examples include periodic, quasi-periodic and random initial data. On the lattice we prove that solutions are polynomially bounded in time for any bounded data. In the continuum, local existence is proved for real analytic data by a Newton iteration scheme.  Global existence for NLS with a regularized nonlinearity follows by analyzing a local energy norm. 

\bigskip(Dedicated to Joel Lebowitz with admiration for his inspiring leadership.)
\end{abstract}
\maketitle

\section{Introduction}
The aim of this note is to study nonlinear Schr{\"o}dinger type equations on $\mathbb{Z}$ and on $\mathbb{R}$ with initial data that are bounded, and do not vanish at infinity. Examples include periodic, quasi-periodic and random data.  We present some modest results describing the dynamics for such data. Although we shall phrase our results in one dimension, most of our methods can be adapted to higher dimensions.\medskip

There is an extensive literature on space periodic data, with bounds on the Sobolev norms as a function of time \cite{MR1215780, MR1309539, MR1777342}.  More recently,  localized perturbations of periodic initial data have been studied \cite{MR3696604, MR3680928}.    The periodic case appears naturally in nonlinear optics, describing periodic signal propagating through a fiber, and local perturbations would correspond to noise and modulation due to information carrying signal. \medskip

Much less is known about quasi-periodic and random cases. The quasi-periodic data arises when there is another periodic signal, not commensurate with the periodic signal. The nonlinearity then will produce a solution with arbitrary number of frequencies. In the case of the KdV equation, Damanik and Goldstein and Binder at al. \cite{MR3486173, MR3859361}  have proved the global existence of uniformly bounded almost periodic solutions  for certain small amplitude quasi-periodic initial data. This remarkable result uses the integrability of KdV  and the fact that the corresponding Schr\"odinger operator has only absolutely continuous spectrum. The case of global existence of large initial data has not yet been understood. For non-integrable NLS, T. Oh \cite{MR3328142}, proved the existence of solutions for quasi-periodic data for short time. He also proved that solutions exist for all time for limit-periodic data, \cite{oh2015global}. The global existence of special space-time quasi-periodic solutions was recently established by W. Wang \cite{wang2019infinite} by using Bourgain's semi-algebraic set methods together with a Newton iteration scheme. It is generally believed that typical solutions to non-integrable equations are not uniformly bounded in time unless the maximum is controlled by a conservation law. Note that in one dimension, the maximum of a solution of NLS  with periodic data is bounded by the energy. \medskip

In this paper we consider the nonlinear Schr{\"o}dinger equation
\begin{equation}\label{1.1}
i \psi_{t} + \partial_{xx} \psi = |\psi|^{2} \psi, \qquad \psi : \mathbb{R}^{1 + 1} \rightarrow \mathbb{C}, \qquad \psi : \mathbb{R} \times \mathbb{Z} \rightarrow \mathbb{C}, \qquad \psi(0,x) = \psi_{0}(x),
\end{equation}
with data which are locally in $H^s$, with uniform bounds over the reals.  The case when $\psi : \mathbb{R} \times \mathbb{Z} \rightarrow \mathbb{C}$ corresponds to a nonlinear Schr{\"o}dinger equation on the one dimensional lattice  with $\partial_{xx} $ given by the discrete Laplacian. The case when $\psi : \mathbb{R} \times \mathbb{R} \rightarrow \mathbb{C}$ corresponds to a nonlinear Schr{\"o}dinger equation in the continuum.
What makes the NLS hard to analyze on $\mathbb{R}$ is the lack of a finite propagation speed. The speed of a signal is proportional to the derivative of the solution. 
In contrast to the nonlinear Schr{\"o}dinger equation, the nonlinear wave equation,
\begin{equation}\label{1.2}
u_{tt} - u_{xx} + u^{3} = 0, \qquad u : \mathbb{R}^{1 + 1} \rightarrow \mathbb{R}, \qquad u(0,x) = u_{0}, \qquad u_{t}(0,x) = u_{1},
\end{equation}
 has a finite propagation speed. This means that at time t the solution $u(t,x)$ depends only on data in the backward light cone $\{ x':|x'-x|\le t\}$. The rest of the data can be set to 0. Hence the solution exists for all time by standard arguments.  The following Proposition shows how finite propagation speed enables one to get bounds on the time evolution in one dimension. 
\begin{proposition}\label{p1.1}
If $u(t,x)$ is a solution to the cubic nonlinear wave equation in one dimension, $(\ref{1.2})$ with initial conditions which are uniformly $C^1(\mathbb{R}) \times C^{0}(\mathbb{R})$, then  
\begin{equation}\label{1.2.1}
|u(t,x)|\le C t^{1/3}.
\end{equation}
\end{proposition}
\noindent \emph{Proof:} For any $x_{0} \in \mathbb{R}$ and $t \in [0, \infty)$, if $\chi(x)$ is a smooth cutoff function,
\begin{equation}\label{1.3}
\chi(x) = \left\{
	\begin{array}{ll}
		1  & \mbox{if } |x| \leq 1 \\
		0 & \mbox{if } |x| > 2,
	\end{array}
\right.
\end{equation}
if $v(t,x)$ is the solution to $(\ref{1.2})$ with $u_{0}$ replaced with $\chi(\frac{x - x_{0}}{T}) u_{0}$ and $u_{1}$ replaced with $\chi(\frac{x - x_{0}}{T}) u_{1}$, then by finite propagation speed, $v(t, x_{0}) = u(t, x_{0})$ for all $0 \leq t \leq T$. The solution to $(\ref{1.2})$ has the conserved energy
\begin{equation}\label{1.4}
E(u, u_{t}) = \frac{1}{2} \int (\partial_{x} u(t,x))^{2} dx + \frac{1}{2} \int (u_{t}(t,x))^{2} dx + \frac{1}{4} \int u(t,x)^{4} dx.
\end{equation}
By direct computation,
\begin{equation}\label{1.5}
E(\chi u_{0}, \chi u_{1}) \lesssim T (\| u_{0} \|_{L^{\infty}}^{2} + \| \partial_{x} u_{0} \|_{L^{\infty}}^{2} + \| u_{1} \|_{L^{\infty}}^{2}).
\end{equation}
Then by the Sobolev embedding theorem,
\begin{equation}\label{1.6}
\aligned
|u(t, x_{0})|^{3} = |v(t, x_{0})|^3 \lesssim  \| \partial_{x} v(t) \|_{L^{2}} \| v(t) \|_{L^{4}}^{2} \lesssim E(v, v_{t}) \lesssim T (\| u_{0} \|_{L^{\infty}}^{2} + \| \partial_{x} u_{0} \|_{L^{\infty}}^{2} + \| u_{1} \|_{L^{\infty}}^{2}).
\endaligned
\end{equation}
This proves the proposition. $\Box$
\medskip

\noindent \textbf{Remark:} If the nonlinearity in $(\ref{1.2})$ is replaced by $u^{2p + 1}$, we would observe a growth rate bounded by $C_{p} t^{\frac{1}{p + 2}}$.\medskip

For the solution to the NLS on the lattice, $\psi : \mathbb{R} \times \mathbb{Z} \rightarrow \mathbb{C}$, it is easy to prove that global solutions exist for all uniformly bounded initial data and that the solution grows at most like $C t^{\frac{1}{2}}$. The nonlinear Schr{\"o}dinger equation on the lattice may be thought of as  having approximate finite propagation speed since derivatives are uniformly bounded. The bound $C t^{\frac{1}{2}}$ may be proved for the nonlinear Schr{\"o}dinger equation on the lattice with either focusing or defocusing nonlinearity. However, observe that the proof of Proposition $\ref{p1.1}$ uses the fact that the wave equation is defocusing. For the defocusing, nonlinear Schr{\"o}dinger equation on the lattice, we may improve the bounds to $C t^{\frac{1}{4}}$. See Sec 2.\medskip

\noindent \textbf{Remark:} The solution to the 1D linear Schr\"odinger with initial data  $\psi_0(x)=\sum_j a_j e^{-(x-j)^2},\, \,j\in \mathbb{Z}$  is given by $$C \sum_j a_j \frac{e^{\frac{-(x-j)^2}{4it+1}}}{(4it+1)^{1/2}}.$$ It is easy to show that if the $|a_{j}|$ are uniformly bounded, then the sup norm of this solution is bounded by $Ct^{1/2}$. The phases in the exponential can be cancelled by the $a_j$ so this is the best one can do. On the other hand, if the $a_j$ are independent complex random variables such that $\mathbb{E}|a_j|^2 = 1$, then $\mathbb{E} |e^{it\Delta}\psi_0|^2(x) \le C$. On the lattice, the upper bound $t^{1/2}$ follows from Proposition 2 and a $t^{1/2}$ lower bound is proved in the Appendix for particular $a_j$.  
  \medskip

We also study a simplified model of NLS in which the nonlinearity is regularized.  The Hamiltonian of the  nonlinear interaction we consider has the form  $ \int |u_{\phi}(x, t)|^{4} dx $ where $u_{\phi}$ denotes the convolution of $u$ with a smooth positive function $\phi(x)$ of compact support. In this case we prove that solutions exist for all time and are polynomially bounded for uniformly smooth initial data.   See Theorem 2 in Sec 3.\medskip

The global existence for NLS with bounded smooth initial data is still open. However in the last section we  extend the local time results of T. Oh for quasi-periodic initial to the real analytic setting by developing a Newton iteration for the short time evolution. This is partly inspired by work of Greene and Jacobowitz \cite{MR283728}  on analytic embedding. In particular we can include data of the form $\sum_j a_j e^{-(x-j)^2}$ where  ${j \in \mathbb Z}$ and $|a_j| \le 1$. The $a_j$ may be random.\medskip

\noindent \textbf{Remark:} Besides global existence a natural question would be to prove that the time average of the energy per unit volume of space is uniformly bounded and then to consider the possibility that the limit of the average exists as time goes to infinity. Partial answer is given in Section two, for the case of the NLS on the lattice.\bigskip

\subsection{Invariant measures for NLS equations.}

In many cases there exists an invariant measure which helps to control and describe the time evolution.
The work of Lebowitz, Rose and Speer \cite{MR939505}  on the Gibbs measure for the focussing NLS on the circle is  the
foundational paper in this field. The construction of invariant measures for the NLS is directly related to the statistical mechanics problem with the NLS energy functional defining the theory. This equilibrium measure may also describe a soliton gas of NLS.  This work was extended by Bourgain \cite{MR1309539, MR1374420} who also studied the time evolution for rough data in the support of this measure.  Recent works of Lebowitz, Mounaix and Wang \cite{lebowitz2013approach} as well as  Carlen, Fr\"ohlich and Lebowitz \cite{carlen2016exponential} describe the rate of equilibration to the Gibbs measure when a suitable noise perturbation is added. More closely related to this note is Bourgain's work on NLS in the defocussing dynamics in a periodic box as the period goes to infinity \cite{MR1777342}. He proves that weak limits  of solutions as  L goes to infinity converge to a unique distributional solution in $C(H^s)$, for $s < \frac{1}{2}$, which depends continuously on initial data in compact space time regions.  \medskip     

Let us now consider the  NLS on a periodic lattice of length L. There is an equilibrium Gibbs measure given by $$Z_L^{-1}\exp(-\sum_j^L |\nabla \phi|^2(j) \pm \sum^L_j |\phi(j)|^4)\,.$$  The +  and - sign are the focussing and defocussing cases respectively.
Let $\left\langle  \cdot \right\rangle$ be the translation invariant expectation and define $\phi(t,j, \omega)$ to be the solution of the periodic NLS with  random initial data  distributed by the Gibbs measure. Then since the total energy is conserved and  $\left\langle |\phi(t, j)|^4  \right\rangle$ is independent of j, we conclude that
$\left\langle |\phi(t, 0)|^4  \right\rangle$ is uniformly bounded independently of L. This shows that averaging gives good control of the time evolution. We speculate that the time average of the local energy at 0 is bounded. In addition one would like to describe its time fluctuations. But this is far beyond the scope of this note.
In \cite{lukkarinen2011weakly},
Lukkarinen and Spohn investigated the dynamics of the lattice NLS with initial data governed by the infinite volume Gibbs measure. They proved that as the coupling goes to 0,  the time rescaled solution to NLS obeys a kinetic equation. Recent work of Mendl and Spohn \cite{mendl2015low}   describes equilibrium time correlations on a one  dimensional lattice. 
On the 3 dimensional lattice,  Chatterjee and Kirkpatrick, \cite{chatterjee2012probabilistic}, studied the statistical mechanics with the focussing non linearity. They prove that as the density is varied, a first order phase transition occurs  corresponding to soliton collapse. \medskip

\noindent \textbf{Remark:} In \cite{turaev2010analytical} the detailed dynamics of a gas of solitons of the Ginzburg-Landau equation is studied. In particular, it is shown that there are solutions in which pairs of coupled solitons and separated from others is a possible class of (chaotic) solutions. The approach in this paper is based on the analysis of the infinite system of coupled ODE's corresponding to the internal degrees of freedom of the soliton pairs (center of mass, relative phase, amplitude). It should be pointed out that this system is not Hamiltonian, as it corresponds to the complex NLS.\medskip







\section{Dynamics of NLS on the lattice }

\medskip The Schr{\"o}dinger equation on the lattice has properties similar to the nonlinear wave equation.
In this section let $-\Delta = \partial^*\partial$ be the finite difference Laplacian on  $\mathbb{Z}$.
Here
\begin{equation}\label{2.1} 
\partial f (x) = f(x+1)-f(x),\, and\,\,\,  \partial^* f(x) = f(x-1)-f(x) 
 \end{equation}

\medskip The lattice NLS is given by
\begin{equation}\label{2.1}
i\frac{\partial}{\partial t}\psi(t,x)= i\dot \psi(t,x) =
-\Delta \psi(t,x) + |\psi|^{2} \psi(t,x), \,\,x \in \mathbb{Z} 
\end{equation}

\begin{proposition}\label{p1}
If $|\psi(0,x)| \le A $ then there exists a constant $C$ such that for any $x_{0} \in \mathbb{Z}$ and $t_{0} \geq 1$, we have $|\psi(t_0, x_{0})|  \le C A \, t_0^{1/2}$, and
\begin{equation}
\frac{1}{t_{0}} \sum_{|x - x_{0}| \leq t_{0}} |\psi(t_{0}, x)|^{2} \lesssim A^{2}.
\end{equation}
\end{proposition}
\noindent \emph{Proof:} The proof of local existence and uniqueness is an application of Picard iteration on the space $l^{\infty}(\mathbb{Z})$ since the linear operator $\Delta$ is bounded on $l^{\infty}(\mathbb{Z})$, the Taylor series for $e^{it \Delta}$ is uniformly bounded for times $0 \leq t \leq 1$.\medskip

To get a bound on the solution at arbitrary time $t_0$ near $x_0$, where $t_{0}$ is possibly large,
we define 
\begin{equation}\label{2.4}
F(t,x)= \frac{[(x-x_0)^2 +1]^{1/2}}{R(2 t_0-t +1 )},
 \end{equation} 
with $R \geq 1$, and let the local mass be given by
\begin{equation}\label{2.5}
M(t)= \sum_x  |\psi(t,x)|^2 e^{-F(x,t)}
\end{equation} 
By using $\partial^*$ is the adjoint of $\partial$, and summation by parts, we have
\begin{equation}\label{2.6} 
d M(t)/dt =  -\sum_x \{|\psi(t,x)|^2 e^{-F(t,x)}( \dot F\,)+ i [\partial \psi \bar \psi - \partial \bar\psi \psi]\,\partial e^{-F(x,t)} \}\
\end{equation}
The first term inside the braces is positive. To estimate the last term note that
\begin{equation}\label{2.7.0}
\partial e^{-F(x,t)} \approx e^{-F(x,t)}\partial F \approx e^{-F(x,t)}/R(2 t_0-t+1)\,.
\end{equation}
Since 
\begin{equation}\label{2.7.1}
|\partial \psi(x) \bar \psi(x) - \partial \bar\psi(x) \psi(x)| \le |\psi(x+1)|^2 +|\psi(x)|^2
\end{equation}
and $$ \sum_x |\psi(x+1)|^2 e^{-F(x,t)} = \sum_x |\psi(x)|^2 e^{-F(x-1,t)}\le \sum_x |\psi(x)|^2 e^{-F(x,t)}(1+\frac{1}{R(2 t_0-t+1)}  )  $$
we get the inequality for $t_0 \ge t$
\begin{equation}\label{2.5}
 d M(t)/dt \le 3 M(t)/R(2t_0-t+1)
\end{equation}   
This implies 
$M(t_0)\leq \frac{3 \ln(2)}{R} M(0),$ and Proposition $\ref{p1}$ follows from

\begin{equation}\label{2.5.1}
\sum_x e^{-\frac{[(x-x_0)^2 +1]^{1/2}}{R(t_{0} + 1)} }|\psi(t_0, x)|^2\le C \sum_x e^{-\frac{[(x-x_0)^2 +1]^{1/2}}{R(2t_0 +1 )} } A^2 \le C R A^2  t_{0}.
\end{equation}
Taking $R = 1$ completes the proof. $\Box$\medskip


If we replace the local mass by the local energy $|\nabla \psi|^2(x) + |\psi(x)|^4$ in (2.5) we can get a improved estimate on the time growth $t^{1/4}$,  compatible with Proposition \ref{p1}. This result only holds for the defocussing case, and would not hold in the focussing case, unlike the proof of Proposition $\ref{p1}$.
\begin{proposition}\label{p2}
If $|\psi(0,x)| \le A $, then for $R \ge 1$,  $|\psi(t_0, x_{0})|  \le C A t_0^{1/4}$, for any $x_{0} \in \mathbb{Z}$, and when $t_{0} \geq 1$,
\begin{equation}
\frac{1}{t_{0}} \sum_{|x - x_{0}| \leq t_{0}} |\psi(t_{0}, x)|^{4} \lesssim A^{4}.
\end{equation}

\end{proposition}
\emph{Proof:} Define the localized energy,
\begin{equation}\label{2.6}
E(t) = \frac{1}{2} \sum_{x} |\psi(t, x + 1) - \psi(t, x)|^{2} e^{-F(t,x)} + \frac{1}{4} \sum_{x} |\psi(t,x)|^{4} e^{-F(t,x)}.
\end{equation}
If we set $\langle v,w \rangle = Re\, v\bar{w}$,\, then by direct computation,
\begin{equation}\label{2.7}
\aligned
\frac{dE}{dt} = -\frac{1}{2} \sum_{x} |\psi(t,x + 1) - \psi(t,x)|^{2} e^{-F(t,x)} \dot{F}(t,x) - \frac{1}{4} \sum_{x} |\psi(t,x)|^{4} e^{-F(t,x)} \dot{F}(t,x) \\
+ \sum_{x} \langle \partial^{\ast} \partial \psi e^{-F}, \dot{\psi} \rangle + \sum_{x} \langle \partial \psi \partial^{\ast} e^{-F}, \dot{\psi} \rangle + \sum_{x} \langle e^{-F} |\psi|^{2} \psi, \dot{\psi} \rangle \\ \leq \sum_{x} \langle i \dot{\psi}, e^{-F} \dot{\psi} \rangle + \sum_{x} \langle \partial \psi \partial^{\ast} e^{-F}, \dot{\psi} \rangle = \sum_{x} \langle \partial \psi \partial^{\ast} e^{-F}, \dot{\psi} \rangle
= \sum_{x} \langle \partial \psi \partial^{\ast} e^{-F}, i \Delta \psi(t,x) - i |\psi|^{2} \psi(t,x) \rangle.
\endaligned
\end{equation}
Using the fact that $\partial$ and $\partial^{\ast}$ are bounded operators on the lattice, $(\ref{2.7.0})$ implies that
\begin{equation}
\sum_{x} \langle \partial \psi \partial^{\ast} e^{-F}, -i |\psi|^{2} \psi \rangle \lesssim \sum_{x} \frac{1}{R(2t_{0} + 1 - t)} e^{-F} |\psi|^{4} \lesssim \frac{1}{R(2t_{0} + 1 - t)} E(t).
\end{equation}
Also, since $\partial$ and $\partial^{\ast}$ are bounded operators,
\begin{equation}
\sum_{x} \langle \partial \psi \partial^{\ast} e^{-F}, i \Delta \psi(t,x) \rangle \lesssim \sum_{x} \frac{1}{R(2t_{0} + 1 - t)} e^{-F(t,x)} |\psi(t, x + 1) - \psi(t,x)|^{2} \lesssim \frac{1}{R(2t_{0} + 1 - t)} E(t).
\end{equation}
In particular, this implies that $E(t_{0}) \lesssim_{R} E(0)$. Therefore,
\begin{equation} \sum_x e^{-\frac{[(x-x_0)^2 +1]^{1/2}}{R(t_{0} + 1)} }|\psi(t_{0}, x)|^{4} \le C \sum_x e^{-\frac{[(x-x_0)^2 +1]^{1/2}}{R(2t_0 +1 )} } A^{4} \le C R A^{4}  t_0\,.  
\end{equation}
Taking $R = 1$ proves the Proposition. $\Box$\medskip

\noindent \textbf{Remark:} In Propositions $\ref{p1}$ and $\ref{p2}$ we can bound the time derivative $|\frac{\partial \psi}{\partial t}|$ by $O(t^{3/2})$ and $O(t^{3/4})$ respectively.\medskip

\noindent \textbf{Remark:} This argument could be generalized to the defocussing equation
\begin{equation}
i\partial/\partial t \,\psi(t,x)= i\dot \psi(t,x) =
-\Delta \psi(t,x) + |\psi|^{p} \psi(t,x), \,\,x \in \mathbb{Z}.
\end{equation}
In that case $|\psi(t, x)| \lesssim t^{\frac{1}{p + 2}}$.

\section{Regularized nonlinearity}
Next, we study the regularized nonlinear Schr{\"o}dinger equation,
\begin{equation}\label{3.1}
i u_{t} + u_{xx} = \mathcal N(u) = \phi \ast (|\phi \ast u|^{2} (\phi \ast u)), \qquad u : \mathbb{R} \times \mathbb{R} \rightarrow \mathbb{C}, \qquad u(0,x) = u_{0},
\end{equation}
where $\phi$ is a real valued, symmetric Schwartz function, and $\ast$ indicates the usual convolution
\begin{equation}\label{3.2}
(f \ast g)(x) = \int f(y) g(x - y) dy.
\end{equation}
\textbf{Remark:} The operator $f \mapsto \phi \ast f$ is a smoothing operator. One important example of such an operator is the Fourier truncation operator. Since the solution to $i u_{t} + u_{xx} = 0$ travels at velocity $\xi$ at frequency $\xi$, truncating the nonlinear term in Fourier space allows us to treat $(\ref{3.1})$ using a finite propagation speed argument.

A solution to $(\ref{3.1})$ conserves the quantities mass,
\begin{equation}\label{3.3}
M(u(t)) = \int |u(t,x)|^{2} dx,
\end{equation}
and energy,
\begin{equation}\label{3.4}
E(u(t)) = \frac{1}{2} \int |\partial_{x} u(t,x)|^{2} dx + \frac{1}{4} \int |\phi \ast u|^{4} dx.
\end{equation}
\noindent \textbf{Remark:} Both quantities could be infinite in this section.

\subsection{Local well-posedness}
Local well-posedness may be proved using perturbative arguments, for data lying in a large function space.
\begin{theorem}\label{t3.1}
If $u_{0}$ lies in $C^{2}(\mathbb{R})$, then $(\ref{1.1})$ has a local solution for $T(\phi, \| u_{0} \|_{L^{\infty}}, \| \nabla^{2} u_{0} \|_{L^{\infty}}) > 0$. The solution $u$ is bounded and uniformly continuous on $[0, T] \times \mathbb{R}$.
\end{theorem}
\emph{Proof:} This proof would work equally well in the focussing or defocussing cases. We first prove
\begin{equation}\label{3.5}
\| e^{it \Delta} u_{0} \|_{L^{\infty}} \lesssim (1 + t^{3/2}) (\| u_{0} \|_{L^{\infty}} + \| \nabla u_{0} \|_{L^{\infty}} + \| \nabla^{2} u_{0} \|_{L^{\infty}}),
\end{equation}
by using stationary phase.
Without loss of generality, it suffices to show that $e^{it \Delta} u_{0}$ is bounded at the origin by the right hand side of $(\ref{3.5})$. Using the stationary phase kernel of $e^{it \Delta}$,
\begin{equation}\label{3.6}
e^{it \Delta} u_{0}(0) = \frac{1}{C t^{1/2}} \int e^{-i \frac{y^{2}}{4t}} u_{0}(y) dy.
\end{equation}
Let $\chi$ be as in $(\ref{1.3})$. Integrating by parts,
\begin{equation}\label{3.7}
\aligned
\frac{1}{C t^{1/2}} \int e^{-i \frac{y^{2}}{4t}} (1 - \chi(y)) u_{0}(y) dy = \frac{1}{C t^{1/2}} \int \frac{2it}{y} \frac{d}{dy} (e^{-i \frac{y^{2}}{4t}}) (1 - \chi(y)) u_{0}(y) dy \\
= O(t^{1/2} \| u_{0} \|_{L^{\infty}}) + C t^{1/2} \int e^{-i \frac{y^{2}}{4t}} \frac{1}{y} (1 - \chi(y)) u_{0}'(y) dy.
\endaligned
\end{equation}
Making another integration by parts argument shows that the second term on the right hand side of $(\ref{3.7})$ is also bounded by the right hand side of $(\ref{3.5})$.\medskip

Now then, by the fundamental theorem of calculus,
\begin{equation}\label{3.8}
\chi(y) u_{0}(y) = \chi(y) u_{0}(0) + \chi(y) (u_{0}(y) - u_{0}(0)) = \chi(y) u_{0}(0) + \chi(y) \int_{0}^{y} u_{0}'(s) ds.
\end{equation}
Since $\chi(y)$ is smooth and compactly supported, $\| \chi(y) u_{0}(0) \|_{H^{1}} \lesssim \| u_{0} \|_{L^{\infty}}$, and therefore by the Sobolev embedding theorem and the fact that $e^{it \Delta}$ is a unitary operator for $L^{2}$-based Sobolev spaces,
\begin{equation}\label{3.9}
\| e^{it \Delta} (\chi(y) u_{0}(0)) \|_{L^{\infty}} \lesssim \| u_{0} \|_{L^{\infty}}.
\end{equation}
Finally, as in $(\ref{3.7})$,
\begin{equation}\label{3.10}
\aligned
\frac{1}{C t^{1/2}} \int e^{-i \frac{y^{2}}{4t}} \chi(y) (u_{0}(y) - u_{0}(0)) dy \\ = \frac{1}{C t^{1/2}} \int \frac{2it}{y} \frac{d}{dy} (e^{-i \frac{y^{2}}{4t}}) \chi(y) (u_{0}(y) - u_{0}(0)) dy \\
= C t^{1/2} \int \frac{d}{dy} (e^{-i \frac{y^{2}}{4t}}) \chi(y) \int_{0}^{1} u_{0}'(sy) ds dy
\endaligned
\end{equation}
The right hand side of $(\ref{3.10})$ is also bounded by the right hand side of $(\ref{3.5})$. This finally proves that $(\ref{3.5})$ holds.\medskip

A similar argument may also be made for the Duhamel term
\begin{equation}\label{3.11}
\int_{0}^{t} e^{i(t - \tau) \Delta} \mathcal N(u) d\tau.
\end{equation}
Since $\phi \ast u$ is a smoothing operator, Theorem $\ref{t3.1}$ also implies
\begin{equation}\label{3.13}
\| e^{i(t - \tau) \Delta} \mathcal N(u) \|_{L^{\infty}} \lesssim_{\phi} (1 + |t - \tau|^{3/2}) \| u \|_{L^{\infty}}^{3}.
\end{equation}

Therefore, Theorem $\ref{t3.1}$ holds by straightforward Picard iteration. To see this, let
\begin{equation}\label{3.14}
X = \{ u : \| u \|_{L_{t,x}^{\infty}([0, T] \times \mathbb{R})} \leq 2 C(\phi) (\| u_{0} \|_{L^{\infty}} + \| \nabla u_{0} \|_{L^{\infty}} + \| \nabla^{2} u_{0} \|_{L^{\infty}}) \},
\end{equation}
where $C(\phi)$ also depends on the implicit constant appearing in $(\ref{3.5})$ when $T \leq 1$. Then define the operator
\begin{equation}\label{3.15}
\Phi(u(t)) = e^{it \Delta} u_{0} - i \int_{0}^{t} e^{i(t - \tau) \Delta} \mathcal N(u) d\tau.
\end{equation}
Then by $(\ref{3.5})$ and $(\ref{3.13})$, if $u \in X$,
\begin{equation}\label{3.16}
\aligned
\| \Phi u(t) \|_{L_{t,x}^{\infty}([0, T] \times \mathbf{R})} \leq C  (\| u_{0} \|_{L^{\infty}} + \| \nabla u_{0} \|_{L^{\infty}} + \| \nabla^{2} u_{0} \|_{L^{\infty}}) \\ + 8 C^{3} T  (\| u_{0} \|_{L^{\infty}} + \| \nabla u_{0} \|_{L^{\infty}} + \| \nabla^{2} u_{0} \|_{L^{\infty}})^{3}.
\endaligned
\end{equation}
Then for $T$ sufficiently small,
\begin{equation}\label{3.17}
\Phi : X \rightarrow X.
\end{equation}
Moreover,
\begin{equation}\label{3.18}
\| \Phi(u(t)) - \Phi(v(t)) \|_{L_{t,x}^{\infty}([0, T] \times \mathbb{R})} \lesssim T \| u - v \|_{L_{t,x}^{\infty}([0, T] \times \mathbb{R})} (\| u \|_{L_{t,x}^{\infty}([0, T] \times \mathbb{R})}^{2} + \| v \|_{L_{t,x}^{\infty}([0, T] \times \mathbb{R})}^{2}),
\end{equation}
which proves that for $T$ sufficiently small, $\Phi$ is a contraction. $\Box$

\subsection{Global well-posedness}
The proof of global well-posedness for $(\ref{3.1})$ relies on mass and energy conservation laws. Although it is not in general true that a function in $L^{\infty}(\mathbb{R})$ will lie in $L^{2}(\mathbb{R})$, it is true that such embeddings do hold on any bounded set. We therefore use local conservation laws to prove global well-posedness.
\begin{theorem}\label{t3.2}
If $u_{0}$ lies in $C^{4}(\mathbb{R})$, then $(\ref{3.1})$ has a global solution whose local energy (see $(\ref{3.19})$ with $R = t^{8}$) is a quantity bounded by $C t^{8}$ when $t$ is large. Moreover, we show that $|u(t,x)| \lesssim t^{8/3}$.
\end{theorem}
\emph{Proof:} Define the local energy
\begin{equation}\label{3.19}
E(x_{0}, t) = \int \chi(\frac{x - x_{0}}{R})^{2} [\frac{1}{2} |u_{x}|^{2} + \frac{1}{4} |\phi \ast u|^{4} + \frac{1}{2} |u|^{2}].
\end{equation}
\noindent \textbf{Remark:} Observe that we rely heavily on the fact that the local energy is positive definite. In the case of finite mass, $(\ref{3.3}) < \infty$, we could prove global well-posedness for both the focussing and defocussing problems. Here, we only prove global well-posedness for the defocussing problem.\medskip

We may compute
\begin{equation}\label{3.20}
\frac{d}{dt} E(x_{0}, t) = \int \chi(\frac{x - x_{0}}{R})^{2} [\langle u_{x}, u_{xt} \rangle + \langle |\phi \ast u|^{2} (\phi \ast u), \phi \ast u_{t} \rangle + \langle u, u_{t} \rangle]
\end{equation}
\begin{equation}\label{3.21}
\aligned
= \int \chi(\frac{x - x_{0}}{R})^{2} [-\langle u_{xx}, u_{t} \rangle + \langle \mathcal N(u), u_{t} \rangle + \langle u, u_{t} \rangle] \\
- \frac{2}{R} \int \chi(\frac{x - x_{0}}{R}) \chi'(\frac{x - x_{0}}{R}) \langle u_{x}, u_{t} \rangle + \int \langle ([\phi, \chi^{2}] u_{t}), |\phi \ast u|^{2} (\phi \ast u) \rangle.
\endaligned
\end{equation}
Here we are using the inner product
\begin{equation}\label{3.22}
\langle f, g \rangle = Re(f(x) \overline{g(x)}),
\end{equation}
and the commutator is given by
\begin{equation}\label{3.23}
[\phi, \chi^{2}] u = \int [\chi^{2}(\frac{y - x_{0}}{R}) - \chi^{2}(\frac{x - x_{0}}{R})] \phi(x - y) u(y) dy.
\end{equation}
Using $(\ref{3.1})$ and the fact that $\langle u_{t}, i u_{t} \rangle = \langle u_{t}, -u_{xx} + \mathcal N(u) \rangle = 0$,
\begin{equation}\label{3.24}
\aligned
\frac{d}{dt} E(x_{0}, t) = \int \chi(\frac{x - x_{0}}{R})^{2} \langle u, u_{t} \rangle - \frac{2}{R} \int \chi(\frac{x - x_{0}}{R}) \chi'(\frac{x - x_{0}}{R}) \langle u_{x}, u_{t} \rangle \\ + \int \langle ([\phi, \chi^{2}] u_{t}), |\phi \ast u|^{2} (\phi \ast u) \rangle = A + B + C.
\endaligned
 \end{equation}\medskip
 
Observe that if $\chi(\frac{x - x_{0}}{R})$ were replaced by $1$, the local energy $(\ref{3.19})$ would be exactly equal to the mass in $\frac{1}{2} \times (\ref{3.3})$ plus the energy in $(\ref{3.4})$. It is also clear that in the case that $\chi = 1$, $\chi'(x) = 0$ and $[\phi, \chi^{2}] = 0$, so the last two terms in $(\ref{3.24})$ would drop out. The first term in $(\ref{3.24})$ would also be zero due to the conservation of mass calculations. For the local energy, we will exploit the fact that a $\frac{1}{R}$ appears in $- \frac{2}{R} \int \chi(\frac{x - x_{0}}{R}) \chi'(\frac{x - x_{0}}{R}) \langle u_{x}, u_{t} \rangle$, and also in $\int \langle ([\phi, \chi^{2}] u_{t}), |\phi \ast u|^{2} (\phi \ast u) \rangle$. Meanwhile, $E(x_{0}, 0) = O(R)$, so terms depending on the size of $u$ itself will grow as $R$ increases. We will show that for any $0 < R < \infty$, we are able to win the tug of war that develops between $\frac{1}{R}$ and the nonlinear terms for a time that increases as $R$ increases. Taking $R \rightarrow \infty$ we shall prove global well-posedness and $|u(t,x)| \lesssim t^{8/3}$.\medskip
   
\noindent \textbf{Term A:} Since $u_{t} = i u_{xx} - i \mathcal N(u)$,
\begin{equation}\label{3.25}
A = \int \chi(\frac{x - x_{0}}{R})^{2} \langle u, u_{t} \rangle = \int \chi(\frac{x - x_{0}}{R})^{2} \langle u, i u_{xx} - i \mathcal N(u) \rangle = A_{1} + A_{2}.
\end{equation}
Integrating by parts,
\begin{equation}\label{3.26}
\aligned
A_{1} = \int \chi(\frac{x - x_{0}}{R})^{2} \langle u, i u_{xx} \rangle = -\frac{2}{R} \int \chi(\frac{x - x_{0}}{R}) \chi'(\frac{x - x_{0}}{R}) \langle u, i u_{x} \rangle \\ \lesssim \frac{1}{R} \| \chi(\frac{x - x_{0}}{R}) u \|_{L^{2}} \| \chi'(\frac{x - x_{0}}{R}) u_{x} \|_{L^{2}} \lesssim \frac{1}{R} E(x_{0}, t)^{1/2} \| \chi'(\frac{x - x_{0}}{R}) u_{x} \|_{L^{2}}.
\endaligned
\end{equation}
By the definition of $\chi$,
\begin{equation}\label{3.27}
|\chi'(x)| \lesssim |\chi(x)| + |\chi(x - 1)| + |\chi(x + 1)|.
\end{equation}
Therefore,
\begin{equation}\label{3.28}
A_{1} \lesssim \frac{1}{R} E(x_{0}, t)^{1/2} \| \chi'(\frac{x - x_{0}}{R}) u_{x} \|_{L^{2}} \lesssim \frac{1}{R} E(x_{0}, t)^{1/2} (\sup_{x_{0}} E(x_{0}, t)^{1/2}).
\end{equation}\medskip

Next,
\begin{equation}\label{3.29}
\aligned
A_{2} = - \int \chi^{2}(\frac{x - x_{0}}{R}) \langle u, i \mathcal N(u) \rangle = \int \langle \phi \ast \chi^{2}(\frac{x - x_{0}}{R})u, i |\phi \ast u|^{2} (\phi \ast u) \rangle
= \int \langle [\phi, \chi^{2}] u, i |\phi \ast u|^{2} (\phi \ast u) \rangle.
\endaligned
\end{equation}
The last equality follows from the fact that
\begin{equation}\label{3.30}
\int \chi^{2}(\frac{x - x_{0}}{R}) Re(-i |\phi \ast u|^{4}) = 0.
\end{equation}

To analyze the commutator $[\phi, \chi^{2}]$, first consider $|x - x_{0}| > 4R$. In that case,
\begin{equation}\label{3.31}
\aligned
 (\ref{3.23}) = [\phi, \chi^{2}] = \int \chi^{2}(\frac{y - x_{0}}{R}) \phi(x - y) u(y) dy \\ = \int \frac{1}{1 + |x - y|^{2}} (1 + |x - y|^{2}) \phi(x - y) \chi^{2}(\frac{y - x_{0}}{R}) u(y) dy.
\endaligned
\end{equation}
Since $\phi$ is a Schwartz function,
\begin{equation}\label{3.32}
\| |x - x_{0}|^{2} \int \chi^{2}(\frac{y - x_{0}}{R}) \phi(x - y) u(y) dy \|_{L^{4}(|x - x_{0}| \geq 4R)} \lesssim_{\phi} \| \chi(\frac{y - x_{0}}{R}) u \|_{L^{2}}.
\end{equation}
Therefore,
\begin{equation}\label{3.33}
\aligned
 \int_{|x - x_{0}| \geq 4R} \langle [\phi, \chi^{2}] u, i |\phi \ast u|^{2} (\phi \ast u) \rangle dx \lesssim_{\phi}  \| \chi(\frac{x - x_{0}}{R}) u \|_{L^{2}} \| \frac{1}{|x - x_{0}|^{2}} |\phi \ast u|^{3} \|_{L^{4/3}(|x - x_{0}| \geq 4R)} \\
 \lesssim E(x_{0}, t)^{1/2} (\sum_{|j| \geq 3} \frac{1}{j^{2} R^{2}} \| \chi(\frac{x - x_{0} - jR}{R}) (\phi \ast u) \|_{L^{4}}^{3}) \lesssim \frac{1}{R} E(x_{0}, t)^{1/2} (\sup_{x_{0}} E(x_{0}, t)^{3/4}).
 \endaligned
\end{equation}

For $|x - x_{0}| \leq 4R$ observe that by the fundamental theorem of calculus,
\begin{equation}\label{3.34}
|\chi^{2}(\frac{y - x_{0}}{R}) - \chi^{2}(\frac{x - x_{0}}{R})| \lesssim \frac{1}{R} |x - y|.
\end{equation}
Then,
\begin{equation}\label{3.35}
\aligned
\int_{|x - x_{0}| \leq 4R} \langle [\phi, \chi^{2}] u, i |\phi \ast u|^{2} (\phi \ast u) \rangle dx \\ \lesssim \frac{1}{R} \| \phi \ast u \|_{L^{4}(|x - x_{0}| \leq 4R)}^{3} \| \int |x - y| |\phi(x - y)| |u(y)| dy \|_{L^{4}(|x - x_{0}| \leq 4R)} \\
\lesssim \frac{1}{R} (\sup_{x_{0}} E(x_{0}, t)^{3/4}) \| \int |x - y| |\phi(x - y)| |u(y)| dy \|_{L^{4}(|x - x_{0}| \leq 4R)}.
\endaligned
\end{equation}
Once again, since $\phi$ is a Schwartz function, by Young's inequality,
\begin{equation}\label{3.36}
\aligned
\| \int |x - y| |\phi(x - y)| |u(y) dy \|_{L^{4}(|x - x_{0}| \leq 4R)} \\ \lesssim \| \int \frac{1}{1 + |x - y|^{2}} |u(y)| dy \|_{L^{4}(|x - x_{0}| \leq 4R)} \lesssim (\sup_{x_{0}} \| \chi(\frac{x - x_{0}}{R}) u \|_{L^{2}}).
\endaligned
\end{equation}
Therefore,
\begin{equation}\label{3.37}
\aligned
 A_2 \le \int \langle [\phi, \chi^{2}] u, i |\phi \ast u|^{2} (\phi \ast u) \rangle \lesssim \frac{1}{R} (\sup_{x_{0}} E(t, x_{0})^{5/4}).
 \endaligned
\end{equation}
\medskip

\noindent \textbf{Term B:} Split
\begin{equation}\label{3.38}
\aligned
- \frac{2}{R} \int \chi(\frac{x - x_{0}}{R}) \chi'(\frac{x - x_{0}}{R}) \langle u_{x}, u_{t} \rangle = - \frac{2}{R} \int \chi(\frac{x - x_{0}}{R}) \chi'(\frac{x - x_{0}}{R}) \langle u_{x}, i u_{xx} \rangle \\ + \frac{2}{R} \int \chi(\frac{x - x_{0}}{R}) \chi'(\frac{x - x_{0}}{R}) \langle u_{x}, i \mathcal N(u) \rangle = B_{1} + B_{2}.
\endaligned
\end{equation}
Decompose
\begin{equation}\label{3.39}
\aligned
B_{1} = - \frac{2}{R} \int \chi(\frac{x - x_{0}}{R}) \chi'(\frac{x - x_{0}}{R}) \langle u_{x}, i u_{xx} \rangle \\ = -\frac{2}{R} \int \chi(\frac{x - x_{0}}{R}) \chi'(\frac{x - x_{0}}{R}) \langle u_{x}, i e^{it \partial_{xx}} u_{xx}(0) \rangle \\ - \frac{2}{R} \int \chi(\frac{x - x_{0}}{R}) \chi'(\frac{x - x_{0}}{R}) \langle u_{x}, i \partial_{xx}(u - e^{it \partial_{xx}} u(0)) \rangle.
\endaligned
\end{equation}

Using $(\ref{3.5})$ and the fact that $\partial_{x}$ commutes with $e^{it \partial_{xx}}$,
\begin{equation}\label{3.40}
\aligned
-\frac{2}{R} \int \chi(\frac{x - x_{0}}{R}) \chi'(\frac{x - x_{0}}{R}) \langle u_{x}, i e^{it \partial_{xx}} u_{xx}(0) \rangle \\ \lesssim \frac{1}{R^{1/2}} \| \chi(\frac{x - x_{0}}{R}) u_{x} \|_{L^{2}} \| e^{it \partial_{xx}} u_{xx}(0) \|_{L^{\infty}} \lesssim (1 + t^{3/2}) \frac{1}{R^{1/2}} E(x_{0}, t)^{1/2} \| u_{0} \|_{C^{4}(\mathbb{R})}.
\endaligned
\end{equation}

Next observe that by Duhamel's formula,
\begin{equation}\label{3.41}
u - e^{it \partial_{xx}} u(0) = i \int_{0}^{t} e^{i(t - \tau) \partial_{xx}} \mathcal N(u) d\tau.
\end{equation}
Next,
\begin{equation}\label{3.42}
\aligned
 \chi'(\frac{x - x_{0}}{R}) \partial_{xx} \int_{0}^{t} e^{i(t - \tau) \partial_{xx}} \mathcal N(u) d\tau = \chi'(\frac{x - x_{0}}{R}) \int_{0}^{t} e^{i(t - \tau) \partial_{xx}} \phi_{xx} \ast (|\phi \ast u|^{2} (\phi \ast u)) d\tau.
 \endaligned
 \end{equation}
 Now then, by Strichartz estimates and the fact that $\phi$ is a Schwartz function, and thus in $L^{1}$,
 \begin{equation}\label{3.43}
 \aligned
\| \chi'(\frac{x - x_{0}}{R}) \int_{0}^{t} e^{i(t - \tau) \partial_{xx}} \phi_{xx} \ast \chi(\frac{x - x_{0}}{4R}) (|\phi \ast u|^{2} (\phi \ast u)) d\tau \|_{L^{2}(\mathbb{R})} \\
\lesssim_{\phi} (\int_{0}^{t} \| \chi(\frac{x - x_{0}}{4R}) (|\phi \ast u|^{2} (\phi \ast u)) \|_{L^{4/3}}^{8/7} d\tau)^{7/8} \lesssim (\int_{0}^{t} (\sup_{x_{0}} E(x_{0}, \tau))^{6/7} d\tau)^{7/8}.
\endaligned
 \end{equation}
 Meanwhile, computing the kernel for a generic function $f$,
\begin{equation}\label{3.44}
\aligned
 \chi'(\frac{x - x_{0}}{R}) \int_{0}^{t} e^{i(t - \tau) \partial_{xx}} \phi_{xx} \ast (1 - \chi(\frac{x - x_{0}}{4R})) f d\tau \\
 =  \chi'(\frac{x - x_{0}}{R}) \int_{0}^{t} \int \int e^{-i \frac{(x - y)^{2}}{4(t - \tau)}} \phi^{(2)}(y - z) (1 - \chi(\frac{z - x_{0}}{4R})) f(z) d\tau dy dz.
 \endaligned
\end{equation}
Since $\chi'(\frac{x - x_{0}}{R})$ is supported on $|x - x_{0}| \leq 2R$ and $(1 - \chi(\frac{z - x_{0}}{4R}))$ is supported on $|z - x_{0}| > 4R$, $|x - z| > 2R$ in $(\ref{3.44})$. In the region where $|y - z| > |x - y|$, using the fact that $\phi$ is a Schwartz function and Young's inequality,
 \begin{equation}\label{3.45}
 \aligned
 \| \chi'(\frac{x - x_{0}}{R}) \int_{0}^{t} \int \int_{|y - z| > |x - y|} e^{-i \frac{(x - y)^{2}}{4(t - \tau)}} \phi^{(2)}(y - z) (1 - \chi(\frac{z - x_{0}}{4R})) f(z) d\tau dy dz \|_{L^{2}} \\
 \lesssim \| \chi'(\frac{x - x_{0}}{R}) \int_{0}^{t} \int \int_{|y - z| > |x - y|} \frac{1}{1 + |x - y|^{2}} \cdot \frac{1}{1 + |x - z|^{2}} \\ \times |\phi^{(2)}(y - z) (y - z)^{4}| (1 - \chi(\frac{z - x_{0}}{4R})) |f(z)| dz dy d\tau \|_{L^{2}} \\
 \lesssim_{\phi} \sum_{j \neq 0} \frac{1}{1 + j^{2} R^{2}} \int_{0}^{t} \| \chi(\frac{x - x_{0} - jR}{R}) f \|_{L^{4/3}} d\tau \lesssim \frac{1}{R^{2}} \int_{0}^{t} (\sup_{x_{0}} \| \chi(\frac{x - x_{0}}{R}) f \|_{L^{4/3}}).
 \endaligned
 \end{equation}

 Next, consider the case when $|x - y| > |y - z|$. In this case we use the fact that
 \begin{equation}\label{3.46}
 e^{-i \frac{(x - y)^{2}}{4(t - \tau)}} = \frac{2i (t - \tau)}{(y - x)} \frac{d}{dy} (e^{-i \frac{(x - y)^{2}}{4(t - \tau)}}).
 \end{equation}
 Then let $\psi$ be a smooth, compactly supported function, $\psi(x) = 1$ for $|x| \leq 1$, and $\psi$ supported on $|x| \leq 2$. Using $(\ref{3.46})$ and integrating by parts,
  \begin{equation}\label{3.47}
 \aligned
  \chi'(\frac{x - x_{0}}{R}) \int_{0}^{t} \int \int \psi(\frac{y - z}{x - y}) e^{-i \frac{(x - y)^{2}}{4(t - \tau)}} \phi^{(2)}(y - z) (1 - \chi(\frac{z - x_{0}}{4R})) f(z) d\tau dy dz \\
=  \chi'(\frac{x - x_{0}}{R}) \int_{0}^{t} \int \int \frac{2i(t - \tau)}{(y - x)^{2}} \psi(\frac{y - z}{x - y}) e^{-i \frac{(x - y)^{2}}{4(t - \tau)}} \phi^{(2)}(y - z) (1 - \chi(\frac{z - x_{0}}{4R})) f(z) d\tau dy dz \\
-  \chi'(\frac{x - x_{0}}{R}) \int_{0}^{t} \int \int \frac{2i(t - \tau)}{(y - x)} \psi'(\frac{y - z}{x - y}) \cdot (\frac{1}{x - y} + \frac{y - z}{(x - y)^{2}}) \\ \times e^{-i \frac{(x - y)^{2}}{4(t - \tau)}} \phi^{(2)}(y - z) (1 - \chi(\frac{z - x_{0}}{4R})) f(z) d\tau dy dz \\
 - \chi'(\frac{x - x_{0}}{R}) \int_{0}^{t} \int \int \frac{2i(t - \tau)}{(y - x)} \psi(\frac{y - z}{x - y}) e^{-i \frac{(x - y)^{2}}{4(t - \tau)}} \phi^{(3)}(y - z) (1 - \chi(\frac{z - x_{0}}{4R})) f(z) d\tau dy dz.
 \endaligned
 \end{equation}
 Applying $(\ref{3.46})$ and integrating by parts three more times, and then making the argument in $(\ref{3.43})$--$(\ref{3.45})$ with $x - y$ and $y - z$ reversed, implies that
 \begin{equation}\label{3.48}
 \| (\ref{3.42}) \|_{L^{2}} \lesssim_{\phi} (\int_{0}^{t} (\sup_{x_{0}} E(x_{0}, \tau))^{6/7} d\tau)^{7/8} \\ +  \frac{(1 + t^{4})}{R^{2}} \int_{0}^{t} (\sup_{x_{0}} \| \chi(\frac{x - x_{0}}{R}) (\phi \ast u) \|_{L^{4}}^{3}) d\tau.
 \end{equation}

Then consider
\begin{equation}\label{3.49}
B_{2} = \frac{2}{R} \int \chi(\frac{x - x_{0}}{R}) \chi'(\frac{x - x_{0}}{R}) \langle u_{x}, i \mathcal N(u) \rangle.
\end{equation}
We can compute
\begin{equation}\label{3.50}
\aligned
B_{2} \lesssim \frac{1}{R} \| \chi(\frac{x - x_{0}}{R}) u_{x} \|_{L^{2}} \| \chi'(\frac{x - x_{0}}{R}) \phi \ast (|\phi \ast u|^{2} (\phi \ast u)) \|_{L^{2}} \\
\lesssim \frac{1}{R} E(t, x_{0})^{1/2} \| \chi'(\frac{x - x_{0}}{R}) \phi \ast (|\phi \ast u|^{2} (\phi \ast u)) \|_{L^{2}}.
\endaligned
\end{equation}
To simplify notation let $f = |\phi \ast u|^{3}$. Since $\phi$ is a Schwartz function,
\begin{equation}\label{3.51}
\aligned
\| \chi'(\frac{x - x_{0}}{R}) \phi \ast f \|_{L^{2}} \lesssim \| \int \chi'(\frac{x - x_{0}}{R}) |\phi(x - y) (x - y)^{2}| \frac{1}{|x - y|^{2}} |f(y)| dy \|_{L^{2}} \\
\lesssim_{\phi} \sum_{j} \frac{1}{1 + j^{2} R^{2}} \| \chi(\frac{x - x_{0} - jR}{R}) f \|_{L^{4/3}} \lesssim (\sup_{x_{0}} \| \chi(\frac{x - x_{0}}{R}) (\phi \ast u) \|_{L^{4}}^{3}).
\endaligned
\end{equation}\medskip

\noindent \textbf{Term C:} For this term as well, split $u_{t} = i u_{xx} - i \mathcal N(u)$.
\begin{equation}\label{3.52}
\aligned
\int \langle ([\phi, \chi^{2}] u_{t}), |\phi \ast u|^{2} (\phi \ast u) \rangle = \int \langle ([\phi, \chi^{2}] i u_{xx}), |\phi \ast u|^{2} (\phi \ast u) \rangle \\ - \int \langle ([\phi, \chi^{2}] i \mathcal N(u)), |\phi \ast u|^{2} (\phi \ast u) \rangle = C_{1} + C_{2}.
\endaligned
\end{equation}
Once again we use the definition of the commutator in $(\ref{3.5})$. Integrating by parts,
\begin{equation}\label{3.53}
\aligned
C_{1} = \int [\chi^{2}(\frac{y - x_{0}}{R}) - \chi^{2}(\frac{x - x_{0}}{R})] \phi(x - y) (i u_{xx}(t,y)) (|\phi \ast u|^{2} (\overline{\phi \ast u}))(t,x) dy dx \\
= -\frac{2}{R} \int \chi(\frac{y - x_{0}}{R}) \chi'(\frac{y - x_{0}}{R}) \phi(x - y) (i u_{x}(t,y)) (|\phi \ast u|^{2} (\overline{\phi \ast u}))(t,x) dy dx \\
- \int [\chi^{2}(\frac{x - x_{0}}{R}) - \chi^{2}(\frac{y - x_{0}}{R})] \phi'(x - y) (i u_{x}(t,y)) (|\phi \ast u|^{2} (\overline{\phi \ast u}))(t,x) dy dx.
\endaligned
\end{equation}
Since $\phi$ is a Schwartz function,
\begin{equation}\label{3.54}
\aligned
\frac{2}{R} \int \chi(\frac{y - x_{0}}{R}) \chi'(\frac{y - x_{0}}{R}) \phi(x - y) (i u_{x}(t,y)) (|\phi \ast u|^{2} (\overline{\phi \ast u}))(t,x) dy dx \\
\lesssim \frac{1}{R} \int |\chi(\frac{y - x_{0}}{R})| |\chi'(\frac{y - x_{0}}{R})| \frac{1}{1 + |x - y|^{2}} \cdot (1 + |x - y|^{2}) |\phi(x - y)| \\ \times |\phi \ast u(t,x)|^{3} |u_{x}(t,y)| dy dx,
\endaligned
\end{equation}
so by Young's inequality,
\begin{equation}\label{3.55}
\lesssim_{\phi} \frac{1}{R} \sum_{j} \frac{1}{1 + j^{2} R^{2}} \| \chi(\frac{x - x_{0} - jR}{R}) |\phi \ast u| \|_{L^{4}}^{3} \| \chi(\frac{x - x_{0}}{R}) u_{x} \|_{L^{2}} \lesssim \frac{1}{R} (\sup_{x_{0}} E(x_{0}, t)^{5/4}).
\end{equation}
Next, consider the term
\begin{equation}\label{3.56}
\int [\chi^{2}(\frac{x - x_{0}}{R}) - \chi^{2}(\frac{y - x_{0}}{R})] \phi'(x - y) (i u_{x}(t,y)) (|\phi \ast u|^{2} (\overline{\phi \ast u}))(t,x) dy dx.
\end{equation}
Integrating by parts again, if the derivative falls on $\chi^{2}(\frac{y - x_{0}}{R})$, then it is possible to proceed as above, only with $u_{x}$ replaced by $u$. Therefore, it is left to estimate
\begin{equation}\label{3.57}
\int [\chi^{2}(\frac{x - x_{0}}{R}) - \chi^{2}(\frac{y - x_{0}}{R})] \phi''(x - y) (i u(t,y)) (|\phi \ast u|^{2} (\overline{\phi \ast u}))(t,x) dy dx.
\end{equation}
To estimate this integral, it is useful to split the integral into three regions, $|x - x_{0}| \leq 4R$ and $|y - x_{0}| \leq 4R$, $|x - x_{0}| > 4R$, and $|y - x_{0}| > 4R$. When $|x - x_{0}| \leq 4R$ and $|y - x_{0}| \leq 4R$, then using $(\ref{3.34})$ and Young's inequality,
\begin{equation}\label{3.58}
\aligned
\int [\chi^{2}(\frac{x - x_{0}}{R}) - \chi^{2}(\frac{y - x_{0}}{R})] \phi''(x - y) (i u(t,y)) (|\phi \ast u|^{2} (\overline{\phi \ast u}))(t,x) dy dx \\
\lesssim_{\phi} \frac{1}{R} \int \int_{|x - x_{0}| \leq 4R, |y - x_{0}| \leq 4R} |x - y| |\phi''(x - y)| |u(t,y)| |\phi \ast u(t,x)|^{3} dx dy \\
\lesssim \sup_{x_{0}} \| \chi(\frac{y - x_{0}}{R}) u \|_{L^{2}} \| \chi(\frac{x - x_{0}}{R}) (\phi \ast u) \|_{L^{4}}^{3} \lesssim \sup_{x_{0}} E(x_{0}, t)^{5/4}.
\endaligned
\end{equation}
If $|x - x_{0}| > 4R$,
\begin{equation}\label{3.59}
\chi^{2}(\frac{x - x_{0}}{R}) - \chi^{2}(\frac{y - x_{0}}{R}) = -\chi^{2}(\frac{y - x_{0}}{R}).
\end{equation}
Then,
\begin{equation}\label{3.60}
\aligned
\int_{|x - x_{0}| > 4R} \chi^{2}(\frac{y - x_{0}}{R}) |\phi''(x - y)| |u(t,y)| |\phi \ast u(t,x)|^{3} dx dy \\
= \int_{|x - x_{0}| > 4R} \chi^{2}(\frac{y - x_{0}}{R}) |\phi''(x - y) (x - y)^{2}| \frac{1}{|x - y|^{2}} |u(t,y)| |\phi \ast u(t,x)|^{3} dx dy \\
\lesssim_{\phi} \sum_{j \neq 0} \frac{1}{j^{2} R^{2}} \| \chi(\frac{x - x_{0} - jR}{R}) u \|_{L^{2}} \| \chi(\frac{y - x_{0}}{R}) |\phi \ast u|^{3} \|_{L^{4/3}} \lesssim \frac{1}{R^{2}} \sup_{x_{0}} E(x_{0}, t)^{5/4}.
\endaligned
\end{equation}
A similar computation may be made for $|y - x_{0}| > 4R$, since in that case,
\begin{equation}\label{3.61}
\chi^{2}(\frac{x - x_{0}}{R}) - \chi^{2}(\frac{y - x_{0}}{R}) = \chi^{2}(\frac{x - x_{0}}{R}).
\end{equation}
Therefore,
\begin{equation}\label{3.62}
C_{1} = \int \langle ([\phi, \chi^{2}] i u_{xx}, |\phi \ast u|^{2} (\phi \ast u) \rangle \lesssim \frac{1}{R} E(x_{0}, t)^{5/4}.
\end{equation}

Finally, consider the contribution of $C_{2}$. If $|y - x_{0}| > 4R$,
\begin{equation}\label{3.63}
\aligned
\int_{|y - x_{0}| > 4R} [\chi^{2}(\frac{y - x_{0}}{R}) - \chi^{2}(\frac{x - x_{0}}{R})] \phi(x - y) \\ \times (\phi \ast |\phi \ast u|^{2} (\phi \ast u)(t,x)) (i |\phi \ast u|^{2} (\overline{\phi \ast u}))(t,y) dy dx \\
= - \int_{|y - x_{0}| > 4R} \chi^{2}(\frac{x - x_{0}}{R}) \tilde{\phi}(x - y) |\phi \ast u|^{2} (\phi \ast u)(t,x) (i |\phi \ast u|^{2} (\overline{\phi \ast u}))(t,y) dy dx
\endaligned
\end{equation}
\begin{equation}\label{3.64}
\aligned
\lesssim \| \chi(\frac{x - x_{0}}{R}) |\phi \ast u|^{2} (\phi \ast u) \|_{L^{4/3}} \| \int_{|y - x_{0}| > 4R} \chi(\frac{x - x_{0}}{R}) |\tilde{\phi}(x - y)| |\phi \ast u|^{3} dy \|_{L^{4}} \\
\lesssim (\sup_{x_{0}} E(x_{0}, t)^{3/4}) \| \int_{|y - x_{0}| > 4R} \chi(\frac{x - x_{0}}{R}) |\tilde{\phi}(x - y)| |\phi \ast u|^{3} dy \|_{L^{4}}.
\endaligned
\end{equation}
Because $\phi$ is a symmetric, real valued, Schwartz function,
\begin{equation}\label{3.65}
\tilde{\phi}(x - y) = \int \phi(x - z) \phi(z - y) dz = \int \phi(x - z) \phi(y - z) dz = \int \phi(x - y - z) \phi(z) dz
\end{equation}
is also a symmetric, real valued, Schwartz function. Therefore,
\begin{equation}\label{3.66}
\aligned
\| \chi(\frac{x - x_{0}}{R}) \int_{|y - x_{0}| > 4R} |\tilde{\phi}(x - y)| |\phi \ast u|^{3} dy \|_{L^{4}} \\ = \| \chi(\frac{x - x_{0}}{R}) \int_{|y - x_{0}| > 4R} |\tilde{\phi}(x - y)| |x - y|^{2} \frac{1}{|x - y|^{2}} |\phi \ast u|^{3} dy \|_{L^{4}} \\
\lesssim_{\phi} \sum_{j \neq 0} \frac{1}{j^{2} R^{2}} \| \chi(\frac{x - x_{0} - jR}{R}) |\phi \ast u|^{3} \|_{L^{4/3}} \lesssim \frac{1}{R^{2}} (\sup_{x_{0}} E(x_{0}, t)^{3/4}).
\endaligned
\end{equation}
By symmetry of arguments, a similar estimate could also be obtained when $|x - x_{0}| > 4R$. Finally, consider the case when $|x - x_{0}| \leq 4R$ and $|y - x_{0}| \leq 4R$. By $(\ref{3.34})$,
\begin{equation}\label{3.67}
\aligned
\int \int_{|x - x_{0}| \leq 4R, |y - x_{0}| \leq 4R}  [\chi^{2}(\frac{y - x_{0}}{R}) - \chi^{2}(\frac{x - x_{0}}{R})] \phi(x - y) \\ \times (\phi \ast |\phi \ast u|^{2} (\phi \ast u)(t,x)) (i |\phi \ast u|^{2} (\overline{\phi \ast u}))(t,y) dy dx \\
\lesssim \frac{1}{R} \int_{|x - x_{0}| \leq 4R} \int_{|y - x_{0}| \leq 4R} |x - y| |\tilde{\phi}(x - y)| |\phi \ast u(t,x)|^{3} |\phi \ast u(t,y)|^{3} dx dy \\
\lesssim_{\phi} \frac{1}{R} (\sup_{x_{0}} E(t, x_{0})^{3/2}).
\endaligned
\end{equation}
Therefore,
\begin{equation}\label{3.68}
C_{2} = \int \langle ([\phi, \chi^{2}] i \phi \ast |\phi \ast u|^{2} (\phi \ast u)), |\phi \ast u|^{2} (\phi \ast u) \rangle \lesssim \frac{1}{R} E(x_{0}, t)^{3/2}.
\end{equation}

Therefore, combining the above estimates on A+B+C we have proved
\begin{equation}\label{3.69}
\aligned
\frac{d}{dt} E(x_{0}, t) \lesssim \frac{1}{R}(\sup_{x_{0}} E(x_{0}, t)) + \frac{1}{R} (\sup_{x_{0}} E(x_{0}, t)^{5/4}) + (1 + t^{3/2}) \frac{E(x_{0}, t)^{1/2}}{R^{1/2}} \| u_{0} \|_{C^{4}} \\
+ (\int_{0}^{t} (\sup_{x_{0}} E(x_{0}, \tau))^{6/7} d\tau)^{7/8} + \frac{1 + t^{4}}{R^{2}} \int_{0}^{t} (\sup_{x_{0}} E(x_{0}, \tau)^{3/4}) d\tau + \frac{1}{R} (\sup_{x_{0}} E(x_{0}, t)^{3/2}).
\endaligned
\end{equation}
By H{\"o}lder's inequality,
\begin{equation}\label{3.70}
\sup_{x_{0}} E(x_{0}, 0) \lesssim R.
\end{equation}
Making a bootstrap argument, suppose $[0, T]$ is an interval for which $\sup_{t,x_{0}} E(x_{0}, t) \lesssim R$, and also $T \leq R^{1/8}$. Then for any $t \in [0, T]$,
\begin{equation}\label{3.71}
(\ref{3.69}) \lesssim 1 + R^{1/4} + (1 + t^{3/2}) \| u_{0} \|_{C^{4}} + R^{3/4} t^{7/8} + \frac{1 + t^{5}}{R^{5/4}} + R^{1/2}.
\end{equation}
Integrating the right hand side of $(\ref{3.71})$ on the interval $[0, T]$,
\begin{equation}\label{3.72}
\int_{0}^{T} 1 + R^{1/4} + (1 + t^{3/2}) \| u_{0} \|_{C^{4}} + R^{3/4} t^{7/8} + \frac{1 + t^{5}}{R^{5/4}} + R^{1/2} dt \lesssim R^{63/64} + R^{5/32} \| u_{0} \|_{C^{4}} \ll R.
\end{equation}
Making a bootstrap argument, this proves $\sup_{x_{0}} E(x_{0}, t) \lesssim R$ for all $t \in [0, T]$. This proves Theorem $\ref{t3.2}$. Using $(\ref{1.6})$, we obtain $|u(t,x)| \lesssim t^{8/3}$. $\Box$\medskip

\noindent \textbf{Remark:} If we assumed smooth initial data, we could combine the previous proof with Gronwall's inequality to prove bounds on higher order derivatives of the solution.

\section{Real analytic local well-posedness}
In this section we prove a local result for the one dimensional, nonlinear Schr{\"o}dinger equation
\begin{equation}\label{4.1}
i \frac{\partial \psi}{\partial t} = - \Delta \psi + |\psi|^{2} \psi, \qquad \psi(0,x) = \psi_{0}.
\end{equation}
T. Oh \cite{MR3328142} proved local well-posedness for almost periodic data. As a simple example consider 
\begin{equation}\label{4.2}
\psi_{0}(x) = \cos(x) + \cos(\sqrt{2} x).
\end{equation}
Such an initial value problem cannot be solved directly using the usual methods of Strichartz estimates combined with the fact that the initial data is in an $L^{2}$-based space. Moreover, the fact that the data $(\ref{4.2})$ consists of two periodic functions whose periods are not rational multiples of one another prevents the study of $(\ref{4.1})$ on a torus.\medskip

In this note, we study $(\ref{4.1})$ for initial data that is real analytic in a strip of width $3$, that is, the set
\begin{equation}\label{4.3}
\Omega = \{ z : z = x + i y, |y| < 3 \}.
\end{equation}
For example, the initial data
\begin{equation}\label{4.4}
\psi_{0}(x) = \sum_{n \in \mathbf{Z}} a_n e^{-(x - n)^{2}}, \,\, |a_n| \le 1
\end{equation}
satisfies this condition. In fact, $(\ref{4.4})$ is an entire function since if $z = x + iy$,
\begin{equation}\label{4.5}
\psi_{0}(z) = \sum_{n \in \mathbf{Z}} a_n e^{-(z - n)^{2}} = e^{y^{2}} \sum_{n \in \mathbf{Z}} a_n e^{-(x - n)^{2}} e^{-2iy (x - n)},
\end{equation}
so the sum $(\ref{4.5})$ converges uniformly in any compact subset of $\mathbb{C}$, proving that $\psi_{0}(x)$ is an entire function.\medskip

We prove the local well-posedness result for real analytic data. We  suspect that Gevrey or very smooth data will suffice to get the existence of local solutions by using Moser's regularization \cite{MR132859}. See H\"ormander
\cite{MR802486} for an abstract version of the Newton-Nash implicit function.  

\begin{theorem}\label{t4.1}
The initial value problem $(\ref{4.1})$ is locally well-posed on some interval $[-T(\psi_{0}), T(\psi_{0})]$, where $\psi_{0}$ is real analytic in the strip $(\ref{4.3})$.
\end{theorem}
\emph{Proof:} Define the formal iteration as follows. First let
\begin{equation}\label{4.6}
\psi_{1}(t,x) = e^{it \Delta} \psi_{0}.
\end{equation}
Next, define $\xi_{2}(t,x)$, to be the linear correction to $\psi_1$ given by
\begin{equation}\label{4.7}
i \frac{\partial \xi_{2}}{\partial t} = -\Delta \xi_{2} + 2 |\psi_{1}|^{2} \xi_{2} + \psi_{1}^{2} \bar{\xi}_{2} + R_{1}, \qquad \xi_{2}(0,x) = 0,
\end{equation}
and
\begin{equation}\label{4.8}
R_{1} = -i \frac{\partial \psi_{1}}{\partial t} - \Delta \psi_{1} + |\psi_{1}|^{2} \psi_{1} = |\psi_{1}|^{2} \psi_{1}.
\end{equation}
and set
\begin{equation}\label{4.9}
\psi_{2}(t,x) = \psi_{1}(t,x) + \xi_{2}(t,x).
\end{equation}
In general, iteratively define
\begin{equation}\label{4.10}
\psi_{n + 1} = \psi_{n} + \xi_{n + 1},
\end{equation}
where
\begin{equation}\label{4.11}
i \frac{\partial \xi_{n + 1}}{\partial t} = - \Delta \xi_{n + 1} + 2 |\psi_{n}|^{2} \xi_{n + 1} + \psi_{n}^{2} \bar{\xi}_{n + 1} + R_{n}, \qquad \xi_{n + 1}(0,x) = 0,
\end{equation}
and when $n > 1$,
\begin{equation}\label{4.12}
R_{n} = -i \frac{\partial \psi_{n}}{\partial t} - \Delta \psi_{n} + |\psi_{n}|^{2} \psi_{n} = 2 |\xi_{n}|^{2} \psi_{n - 1} + \xi_{n}^{2} \bar{\psi}_{n - 1} + |\xi_{n}|^{2} \xi_{n}.
\end{equation} The rapid convergence of the Newton iteration relies on the fact that $R_n$ is quadratic in $\xi_{n}$.
For any $n > 1$, arguing by induction,
\begin{equation}\label{4.13}
i \frac{\partial \psi_{n}}{\partial t} + \Delta \psi_{n} = i \frac{\partial}{\partial t} (\psi_{1} + \sum_{m = 2}^{n} \xi_{m}) + \Delta (\psi_{1} + \sum_{m = 2}^{n} \xi_{m}) = |\psi_{n}|^{2} \psi_{n} - R_{n + 1}.
\end{equation}
Therefore, if we can show that $\xi_{n} \rightarrow 0$ and $R_{n} \rightarrow 0$ in a suitable Banach space as $n \rightarrow \infty$, we are done.\medskip

It is convenient to rewrite $(\ref{4.11})$ and $(\ref{4.12})$ in matrix notation. Let
\begin{equation}\label{4.14}
u_{n} = \begin{pmatrix} 
\xi_{n} \\ \bar{\xi}_{n}
\end{pmatrix}, \,\,
b_{n} = \begin{pmatrix}
R_{n} \\ -\bar{R}_{n}
\end{pmatrix},\,\,
M_{0} = \begin{pmatrix}
-\Delta & 0 \\
0 & \Delta
\end{pmatrix},\,\,
V_{n}(t,x) = 
\begin{pmatrix}
2 |\psi_{n}|^{2} &  \psi_{n}^{2}  \\
-\bar{\psi}_{n}^{2}  &  -2 |\psi_{n}|^{2} \\
\end{pmatrix}.
\end{equation}
Then $(\ref{4.11})$ has the form
\begin{equation}\label{4.15}
i \frac{\partial u_{n + 1}}{\partial t} = M_{n}(t,x) u_{n + 1} + b_{n}, \qquad M_{n} = M_{0} + V_{n}.
\end{equation}
Let $\Phi(t, t_{1})$ be the fundamental solution operator to this equation,
\begin{equation}\label{4.16}
i \frac{\partial \Phi(t, t_{1})}{\partial t} = M_{n} \Phi_{n}(t, t_{1}), \qquad \Phi_{n}(t_1, t_{1}) = Id.
\end{equation}

Formally, $\Phi_{n}(t, t_{1})$ can be expressed as a time ordered integral
\begin{equation}\label{4.17}
\Phi_{n}(t, t_{1}) = \mathbf{T} \exp(-i \int_{t_{1}}^{t} M(s) ds).
\end{equation}
Then the solution to $(\ref{4.15})$ is given by
\begin{equation}\label{4.18}
u_{n + 1} = -i \Phi_{n}(t, 0) \int_{0}^{t} \Phi_{n}(s, 0)^{-1} b_{n}(s) ds.
\end{equation}
Indeed, under $(\ref{4.18})$, $u_{n + 1}(0) = 0$, and
\begin{equation}\label{4.19}
i \frac{\partial u_{n + 1}}{\partial t} = b_{n}(t) + M_{n} u_{n + 1}.
\end{equation}
\begin{lemma}[Properties of $\Phi$]\label{l4.1}
Let 
\begin{equation}\label{4.20}
|V_{n}(t)| = \sup_{s \leq t, x} |V(s,x)|.
\end{equation}
Then,
\begin{equation}\label{4.21}
\| \Phi_{n} f \|_{L^{2}(\mathbb{R})} \leq e^{|t| |V_{n}(t)|} \| f \|_{L^{2}(\mathbb{R})}.
\end{equation}
\end{lemma}
\emph{Proof:} $\Phi_{n}$ is the solution operator to $(\ref{4.15})$ with $b_{n} = 0$. Standard local well-posedness arguments show that $(\ref{4.15})$ has a solution. Next, $(\ref{4.21})$ holds for the interval of existence via Gronwall's inequality. The estimates $(\ref{4.20})$ and $(\ref{4.21})$ combine to show that $\Phi_{n}(t, t_{1})$ is well-defined for all $t, t_{1} \in \mathbb{R}$. $\Box$

\begin{lemma}[More properties of $\Phi$]\label{l4.2}
\begin{equation}\label{4.22}
[\frac{d}{dx}, \Phi_{n}(t_{2}, t_{1})] = -i \int_{t_{1}}^{t_{2}} \Phi_{n}(t_{2}, s) \frac{d}{dx} V_{n}(s, x) \Phi_{n}(s, t_{1}) ds,
\end{equation}
and
\begin{equation}\label{4.23}
[x, \Phi_{n}(t_{2}, t_{1})] = 2i \int_{t_{1}}^{t_{2}} \Phi_{n}(t_{2}, s) \frac{d}{dx} \Phi_{n}(s, t_{1}) ds.
\end{equation}
\end{lemma}
\emph{Proof:} Without loss of generality suppose $t_{1} = 0$ and $t_{2} = 1$. Then for any $N$,
\begin{equation}\label{4.24}
\Phi_{n}(1, 0) = \Phi_{n}(1, 1 - \frac{1}{N}) \circ \cdots \circ \Phi_{n}(\frac{1}{N}, 0).
\end{equation}
Therefore, by direct computation,
\begin{equation}\label{4.25}
[\frac{d}{dx}, \Phi_{n}(1, 0)] = \sum_{j = 1}^{N} \Phi_{n}(1, 1 - \frac{1}{N}) \circ \cdots \circ [\frac{d}{dx}, \Phi_{n}(\frac{j}{N}, \frac{j - 1}{N})] \circ \cdots \circ \Phi_{n}(\frac{1}{N}, 0).
\end{equation}
Then compute,
\begin{equation}\label{4.26}
[\frac{d}{dx}, \Phi_{n}(\frac{j}{N}, \frac{j - 1}{N})] \sim -\frac{i}{N} \frac{d}{dx} V_{n}(\frac{j - 1}{N}, x).
\end{equation}
Similarly,
\begin{equation}\label{4.27}
[x, \Phi_{n}(\frac{j}{N}, \frac{j - 1}{N})] \sim 2i \frac{1}{N} \frac{d}{dx}.
\end{equation}
This proves the lemma. $\Box$\medskip

By taking the time derivative of $\Phi(t)^{-1} \Phi(t) = Id$,
\begin{equation}\label{4.28}
\frac{d}{dt} (\Phi(t)^{-1}) = -i \Phi(t)^{-2} \Phi(t) M(t) = -i \Phi(t)^{-1} M(t), \qquad \Phi(t) = \mathbf{T} \exp(-i \int_{0}^{t} M(s) ds).
\end{equation}
Next, if $f$ is a real analytic function, making a Taylor expansion,
\begin{equation}\label{4.29}
e^{a \frac{d}{dx}} f = \sum_{n = 0}^{\infty} \frac{a^{n}}{n!} f^{(n)}(x) = f(x + a).
\end{equation}
Therefore, we have the identity,
\begin{equation}\label{4.30}
e^{a \frac{d}{dx}} \Phi_{n, 0} f = \Phi_{n, a} e^{a \frac{d}{dx}} f,
\end{equation}
where $\Phi_{n, a}$ is the fundamental solution for $M_{0} + V_{n}(t, x + a)$.\medskip

The norm that will be used to control the iteration is given by
\begin{equation}\label{4.31}
\| f \|(r, p) = \sup_{x} \sum_{n \geq 0, 0 \leq q \leq p} |\frac{f^{(n + q)}(x)}{n!}| r^{n}.
\end{equation}
The following general lemma relates different norms.
\begin{lemma}\label{l4.3}
Let $r_{n + 1} < r_{n}$ and set $\delta_{n} = r_{n} - r_{n + 1}$. Then,
\begin{equation}\label{4.32}
\| f \|(r_{n + 1}, p) \lesssim_{p} \| f \|(r_{n}, 0) \delta_{n}^{-p}.
\end{equation}
\end{lemma}
\emph{Proof:} We consider the simplest form of the estimate
\begin{equation}\label{4.33}
\| f \|(r_{n + 1}, 1) \lesssim \| f \|(r_{n}, 0) \delta_{n}^{-1}.
\end{equation}
 
\begin{equation}\label{4.34}
\| f \|(r_{n + 1}, 1)= \sup_{x} \sum_{m \geq 0} \frac{|f^{(m + 1)}(x)|}{(m + 1)!} r_{n}^{m + 1} [(m + 1) \frac{r_{n + 1}^{m}}{r_{n}^{m+1}}] .
\end{equation}
The inequality follows since the factor in brackets is uniformly bounded by $ \delta_{n}^{-1}$. The case for general p may proved the same way.  
$\Box$\medskip

\noindent \textbf{Remark:} It may be possible to refine the estimates to bound the implicit constant in $(\ref{4.32})$ to $q!$, which would agree with Cauchy integral formula. However, it is not too important to determine the exact constant here.\medskip

The proof of the convergence of the Newton iteration relies on the bound for the operator $\Phi_{n}(t, s)$, where $|t|, |s| \leq 1$.
\begin{lemma}\label{l4.4}
For $|t| \leq 1$ and any $r_{n + 1} > 0$,
\begin{equation}\label{4.37}
\| \Phi_{n}(t, 0) b_{n} \|(r_{n + 1}, 0) \lesssim e^{t |V_{n}(t)|} \| b_{n} \|(r_{n + 1}, 3)(1 + \| V_{n} \|(r_{n + 1}, 3) + \| V_{n} \|(r_{n + 1}, 3)^{2}).
\end{equation}
\end{lemma}
\emph{Proof:} As in the previous section we have
\begin{equation}\label{4.38}
\| e^{it \Delta} \psi_{0} \|(r,0) \lesssim (1 + t^{3/2}) (\| \psi_{0} \|(r,0) + \| \nabla \psi_{0} \|(r,0) + \| \nabla^{2} \psi_{0} \|(r,0).
\end{equation}
This takes care of the case when $V_{n} \equiv 0$.
When $V_{n}$ is not identically equal to zero, we use $(\ref{4.30})$ and the commutator estimates in Lemma $\ref{l4.2}$. Let $\chi_{j}$ be a partition of unity, where $\sum_{j} \chi(x - j) = \sum_{j} \chi_{j}(x) = 1$ for any $x \in \mathbf{R}$, and set
\begin{equation}\label{4.39}
g_{n, j} = [1 + (x - j)^{2}] \Phi_{n}(t, 0) \chi_{j} b_{n}.
\end{equation}
Therefore,
\begin{equation}\label{4.40}
\| \Phi_{n}(t, 0) b_{n} \|(r_{n}, 0) = \| \sum_{j} (1 + (x - j))^{-2} g_{n, j} \|(r_{n}, 0) \leq \sup_{j} \| g_{n, j} \|(r_{n}, 0).
\end{equation}
To simplify notation, let $f = \chi_{j} b_{n}$. Also, without loss of generality suppose that $j = 0$. Then,
\begin{equation}\label{4.41}
\aligned
x^{2} \Phi_{n}(t, 0) f = x \Phi_{n}(t, 0) xf + x \int_{0}^{t} \Phi_{n}(t, s) \frac{d}{dx} \Phi_{n}(s, 0) f ds \\
= x \Phi_{n}(t, 0) xf + x \int_{0}^{t} \Phi_{n}(t, 0) (\frac{d}{dx} f) ds \\ + x \int_{0}^{t} \Phi_{n}(t, s) \int_{0}^{s} \Phi_{n}(s, s_{1}) (\frac{d}{dx} V_{n}(s_{1}, s)) \Phi_{n}(s_{1}, 0) f ds_{1} ds \\
= x \Phi_{n}(t, 0) xf + tx \Phi_{n}(t, 0) (\frac{d}{dx} f) \\ + x \int_{0}^{t} \int_{0}^{s} \Phi_{n}(t, s_{1}) (\frac{d}{dx} V_{n}(s_{1}, s)) \Phi_{n}(s_{1}, 0) f ds_{1} ds.
\endaligned
\end{equation}

Making a similar calculation,
\begin{equation}\label{4.42}
\aligned
x \Phi_{n}(t, 0) xf = \Phi_{n}(t, 0) x^{2} f + t \Phi_{n}(t, 0) (\frac{d}{dx} (xf)) \\ + \int_{0}^{t} \int_{0}^{s} \Phi_{n}(t, s_{1}) (\frac{d}{dx} V_{n}(s_{1}, s)) \Phi_{n}(s_{1}, 0) (xf) ds_{1} ds.
\endaligned
\end{equation}
Then, remembering that $f = \chi_{j} b_{n}$, since $0 \leq t \leq 1$,
\begin{equation}\label{4.43}
\aligned
\| x \Phi_{n}(t, 0) xf \|_{L^{2}} \leq e^{t |V_{n}(t)|} \| x^{2} f \|_{L^{2}} + t e^{t |V_{n}(t)|} \| f \|_{L^{2}} + \\ t e^{t |V_{n}(t)|} \| x \frac{d}{dx} f \|_{L^{2}} + e^{t |V_{n}|} |\frac{d}{dx} V_{n}| \| x f \|_{L^{2}} \lesssim e^{|V_{n}(1)|} (1 + |\frac{d}{dx} V_{n}(1)|) \| b_{n} \|(r_{n + 1}, 1).
\endaligned
\end{equation}
By a similar calculation, since $0 \leq t \leq 1$,
\begin{equation}\label{4.44}
\| t \Phi_{n}(t, 0) (\frac{d}{dx} (xf)) \|_{L^{2}} \lesssim e^{|V_{n}(1)|} (1 + |\frac{d}{dx} V_{n}(1)|) \| b_{n} \|(r_{n + 1}, 2).
\end{equation}

Finally,
\begin{equation}\label{4.45}
\aligned
x \int_{0}^{t} \int_{0}^{s} \Phi_{n}(t, s_{1}) (\frac{d}{dx} V_{n}(s_{1}, s)) \Phi_{n}(s_{1}, 0) f ds_{1} ds \\
= \int_{0}^{t} \int_{0}^{s} \Phi_{n}(t, s_{1}) (\frac{d}{dx} V_{n}(s_{1}, s)) x \Phi_{n}(s_{1}, 0) f ds_{1} ds \\ + \int_{0}^{t} \int_{0}^{s} \int_{s}^{t} \Phi_{n}(t, \tau) \frac{d}{dx} \Phi_{n}(\tau, s_{1}) (\frac{d}{dx} V_{n}(s_{1}, s)) \Phi_{n}(s_{1}, 0) f ds_{1} ds.
\endaligned
 \end{equation}
Also, by $(\ref{4.22})$,
\begin{equation}\label{4.46}
\aligned
\int_{0}^{t} \int_{0}^{s} \int_{s}^{t} \Phi_{n}(t, \tau) \frac{d}{dx} \Phi_{n}(\tau, s_{1}) (\frac{d}{dx} V_{n}(s_{1}, s)) \Phi_{n}(s_{1}, 0) f ds_{1} ds \\
= \int_{0}^{t} \int_{0}^{s} \int_{s}^{t} \Phi_{n}(t, \tau) \Phi_{n}(\tau, s_{1}) \frac{d}{dx} [(\frac{d}{dx} V_{n}(s_{1}, s)) \Phi_{n}(s_{1}, 0) f] ds_{1} ds d\tau \\
+ \int_{0}^{t} \int_{0}^{s} \int_{s}^{t} \Phi_{n}(t, \tau) \int_{s_{1}}^{\tau} \Phi_{n}(\tau, \tilde{\tau}) (\frac{d}{dx} V_{n}(\tilde{\tau}, x)) \Phi_{n}(\tilde{\tau}, s_{1}) (\frac{d}{dx} V_{n}(s_{1}, s)) \Phi_{n}(s_{1}, 0) f ds_{1} ds d\tau d\tilde{\tau}.
\endaligned
\end{equation}
The $L^{2}$ norm of the last term is bounded by
\begin{equation}\label{4.47}
e^{t |V_{n}(t)|} \| V_{n} \|(r_{n + 1}, 3)^{2} \| f \|_{L^{2}}.
\end{equation}
The other remaining terms,
\begin{equation}\label{4.48}
\int_{0}^{t} \int_{0}^{s} \Phi_{n}(t, s_{1}) (\frac{d}{dx} V_{n}(s_{1}, s)) x \Phi_{n}(s_{1}, 0) f ds_{1} ds,
\end{equation}
and
\begin{equation}\label{4.49}
\int_{0}^{t} \int_{0}^{s} \int_{s}^{t} \Phi_{n}(t, \tau) \Phi_{n}(\tau, s_{1}) \frac{d}{dx} [(\frac{d}{dx} V_{n}(s_{1}, s)) \Phi_{n}(s_{1}, 0) f] ds_{1} ds d\tau,
\end{equation}
may be handled in a similar manner.

Finally, using the Sobolev embedding theorem,
\begin{equation}\label{4.50}
\| g_{n,j} \|_{L^{\infty}} \lesssim \| g_{n,j} \|_{L^{2}} + \| \nabla g_{n,j} \|_{L^{2}},
\end{equation}
completes the proof of the lemma. $\Box$\medskip

\emph{Newton iteration:} The induction argument makes use of the fact that if $b_{n}$ is real analytic, then $u_{n + 1}$, given by $(\ref{4.18})$, is also real analytic. Additionally, a product of two real analytic functions is real analytic, as well as the complex conjugate of a real analytic function is real analytic.\medskip

We start the induction by assuming
\begin{equation}\label{4.51}
\epsilon_{1} = \| u_{1} \|(r_{1}, 0) \ll 1,
\end{equation}
and define
\begin{equation}\label{4.52}
\epsilon_{k} = \| u_{k} \|(r_{k}, 0), \qquad k \geq 2.
\end{equation}
It is always possible to assume $(\ref{4.51})$ after rescaling. By Lemmas $\ref{l4.3}$ and $\ref{l4.4}$,
\begin{equation}\label{4.53}
\epsilon_{2} = \| u_{2} \|(r_{2}, 0) \leq C \| \psi_{1} \|(r_{2}, p)^{3} \leq C[\| \psi_{1} \|(r_{1}, 0) \delta_{1}^{-p}]^{3} = [C \epsilon_{1} \delta_{1}^{-p}]^{3}.
\end{equation}
In general, the remainder $b_{n}$ is quadratic in $u_{n}$ by $(\ref{4.13})$ and therefore, for any $n$,
\begin{equation}\label{4.54}
\epsilon_{n + 1} \leq C [C \epsilon_{n} \delta^{-p}]^{2} (C \epsilon_{1} \delta^{-p}).
\end{equation}
If $\epsilon_{1} > 0$ is sufficiently small, $p=3$ and $\delta_{n} = n^{-2}$,
\begin{equation}\label{4.55}
\epsilon_{n + 1} \leq C \epsilon_{n}^{2} n^{4} \leq C^{n} (n!)^{4} \epsilon_{1}^{2^{n}} \rightarrow 0,
\end{equation}
as $n \rightarrow \infty$. Thus, $\psi_{n}$ converges as $n \rightarrow \infty$.

\section{Appendix: Linear Schr\"odinger time evolution on $\mathbb{Z}$ } 

Let $ \psi(t,n) = e^{it \Delta}\psi_0(n)$ where $\psi_0 = \sum_j a_j \delta_j$ and $\Delta$ is the finite difference Laplacian.
\medskip

{\bf Lemma A}  We can choose $a_n,\,|a_j| \le 1$  so that $|\psi(t,0)| \ge \delta t^{1/2}, \delta >0$. \medskip

{\it Proof:} The fundamental solution to the lattice Schr\"odinger equation can be expressed in terms of the integral
$$F_n(t) = (2\pi)^{-1}\int_0^{2\pi}e^{it\cos (\theta)+i n\theta } d\theta, \,\, n\in \mathbb{Z} \,,$$ which is closely related to the Bessel function. This integral  has two saddle points $\theta_s, \, \pi - \theta_s$ where $\sin(\theta_s) = n/t,\,\, \cos(\theta_s) = [1-(n/t)^2]^{1/2}$. By classical stationary phase, if $|n| \le t/2$
 $$F_n(t) \approx [t\cos(\theta_s)]^{-1/2}\, \cos(\phi(t,n)(1 +O(1/t))$$ when n is even and  $$F(t,n)=i[t\cos(\theta_s)]^{-1/2}\sin(\phi(t,n)(1+O(1/t)),$$ when n is odd. Here $$\phi(t,n)= \pi/4+t\cos(\theta_s)+n \theta_s\,.$$
Since  $|\psi(t,0)|= |\sum_n a_n F_n(t)|$ we will choose $a_n$ so that the sum equals
 $$\sum_{ |n|\le t/2}|F(t,n)|  \ge \delta t^{1/2} \,.$$  
To prove this lower bound note that if $|\cos(\phi(t,n))|$ is small when n is even then $|\sin(\phi(t,n+1))|\ge 1/4$  since $|\phi(t,n+1)- \phi(t,n)| \approx |\theta_s|= \arcsin n/t \le \pi/6$.\medskip


\section*{Acknowledgements}
The authors  thank J. Bourgain, P. Deift, J. Lebowitz and W. Schlag for helpful discussions. The first author gratefully acknowledges the support of NSF grants DMS-1500424 and DMS-1764358 while writing this paper. He also gratefully acknowledges the support by the von Neumann fellowship at the Institute for Advanced Study while writing this paper. The second author is supported in part by NSF grant DMS-160074.

\bibliography{biblio}
\bibliographystyle{alpha}
\medskip



\end{document}